\long\def\symbolfootnote[#1]#2{\begingroup\def\thefootnote{\fnsymbol{footnote}}\footnote[#1]{#2}\endgroup}

\documentclass{amsart}
\usepackage{amssymb}
\usepackage{amsmath}
\usepackage{amsfonts}

\newtheorem{thm}{Theorem}[section]
\newtheorem{cor}[thm]{Corollary}

\newtheorem{prop}[thm]{Proposition}

\theoremstyle{remark}
\newtheorem{remark}[thm]{Remark}
\numberwithin{equation}{section}

\theoremstyle{definition}

\begin{document}
\author{Ovidiu Munteanu and Natasa Sesum}
\title{On gradient Ricci solitons}
\date{}

\begin{abstract}
In the first part of the paper we derive integral curvature estimates for
complete gradient shrinking Ricci solitons. Our results and the recent work
in \cite{FG} classify complete gradient shrinking Ricci solitons with
harmonic Weyl tensor. In the second part of the paper we address the issue
of existence of harmonic functions on gradient shrinking K\"{a}hler and
gradient steady Ricci solitons. Consequences to the structure of shrinking
and steady solitons at infinity are also discussed.
\end{abstract}

\maketitle

\section{Introduction and the results}

\symbolfootnote[0]{The first author was partially supported by NSF grant DMS-1005484 and the second author was partially supported by NSF grant DMS-1110145}

A complete Riemannian metric $g$ on a smooth manifold $M$ is called a
gradient Ricci soliton if there is a function $f$ so that 
\begin{equation*}
Ric+Hess\left( f\right) =\rho \cdot g,
\end{equation*}
where $\rho \in \mathbb{R}$. After rescaling the metric $g$ we may assume
that $\rho \in \left\{ -\frac{1}{2},0,\frac{1}{2}\right\} $. Gradient Ricci
solitons arise often as singularity models of the Ricci flow and that is why
understanding them is an important question in the field. Depending on the
behavior of the Ricci flow on solitons, they are called shrinking if $\rho =%
\frac{1}{2},$ steady if $\rho =0$ and expanding if $\rho =-\frac{1}{2}.$

The classification of gradient shrinking Ricci solitons has been a subject
of interest for many people. Hamilton (\cite{Ha1}) showed that the only
closed gradient shrinking Ricci solitons in two dimensions are Einstein. In
three dimensions, Ivey proved that all compact, gradient shrinking Ricci
solitons must have constant positive curvature. The recent work of B\"{o}hm
and Wilking (\cite{BW}) implies the compact gradient shrinking Ricci
solitons with positive curvature operator in any dimension have to be of
constant curvature, generalizing Ivey's result. In higher dimensions, Koiso,
Cao, Feldman, Ilmanen and Knopf constructed examples of gradient shrinking
Ricci solitons that are not Einstein, \cite{C2, FIK}.

The Hamilton-Ivey estimate shows that three dimensional complete solitons
have nonnegative sectional curvatures. Combining this with the results of
Perelman yields that the three dimensional gradient shrinking solitons with
bounded sectional curvatures are $\mathbb{S}^{3}$, $\mathbb{R}^{3}$, $%
\mathbb{S}^{2}\times \mathbb{R}$ and their quotients.

Recently, Ni and Wallach (\cite{NW}) have studied the classification of
complete gradient shrinking Ricci solitons with vanishing Weyl curvature
tensor, in any dimension, under the assumptions of nonnegative Ricci
curvature and at most exponential growth of the norm of curvature operator.
They showed that the only shrinkers satisfying these assumptions are $%
\mathbb{S}^{n}$, $\mathbb{R}^{n}$, $\mathbb{S}^{n-1}\times \mathbb{R}$, and
their quotients. In \cite{CWZ} the assumption on nonnegative Ricci curvature
has been relaxed to having the Ricci curvature bounded from below. Using a
technique developed in the compact setting by \cite{ENM}, Petersen and Wylie
have obtained \cite{PW} the same classification of complete locally
conformally flat gradient shrinking solitons assuming only an integral bound
of the Ricci curvature: 
\begin{equation}
\int_{M}|Ric|^{2}e^{-f}<\infty .  \label{eq-int-ricci}
\end{equation}

In \cite{Z} Zhang proved that gradient shrinking Ricci solitons with
vanishing Weyl tensor must have nonnegative curvature operator, which by any
of the results mentioned above proved the classification of such solitons as
finite quotients of $\mathbb{R}^{n}$, $\mathbb{S}^{n-1}\times \mathbb{R}$ or 
$\mathbb{S}^{n}$.

The question whether certain integral curvature estimates such as (\ref%
{eq-int-ricci}) are true for complete gradient shrinking Ricci solitons has
been raised for example in \cite{PW, CWZ, C}. Besides being interesting on
their own, such estimates would have as a consequence an alternate, simpler
proof of the classification proved in \cite{Z} and should be useful in
proving more general results. In this paper we prove that (\ref{eq-int-ricci}%
) is true for any gradient shrinker. In fact, we establish the following.

\begin{thm}
\label{prop-int-ric0} For any complete gradient shrinking Ricci soliton $%
\left( M,g\right) $ we have 
\begin{equation*}
\int_{M}\left\vert Ric\right\vert ^{2}e^{-\lambda f}<\infty ,\ \ \text{for
any }\lambda >0.
\end{equation*}
\end{thm}

Another curvature quantity which is of interest for classification of
shrinking solitons is $\int_{M}\left\vert \nabla Ric\right\vert ^{2}e^{-f}.$
If this integral is finite, it implies a useful identity (\cite{Cao})%
\begin{equation}
\int_{M}\left\vert \nabla Ric\right\vert ^{2}e^{-f}=\int_{M}\left\vert
div\left( Rm\right) \right\vert ^{2}e^{-f}<\infty ,  \label{iden-curv}
\end{equation}%
which is crucial in the classification result of \cite{CWZ}, mentioned
above. At this time, we do not know if (\ref{iden-curv}) should hold true
for any gradient shrinker. Our next result says that the identity is true
assuming a weighted $L^{2}$ bound of the Riemann curvature tensor. In view
of the right hand side of (\ref{iden-curv}), such an assumption is quite
natural.

\begin{thm}
\label{prop-iden0} Let $(M,g)$ be a gradient shrinking Ricci soliton. If for
some $\lambda <1$ we have $\int_{M}\left\vert Rm\right\vert ^{2}e^{-\lambda
f}<\infty ,$ then the following identity holds: 
\begin{equation*}
\int_{M}\left\vert \nabla Ric\right\vert ^{2}e^{-f}=\int_{M}\left\vert
div\left( Rm\right) \right\vert ^{2}e^{-f}<\infty .
\end{equation*}
\end{thm}

As a consequence, we can prove that (\ref{iden-curv}) is true for gradient
shrinking Ricci solitons with harmonic Weyl tensor. Furthermore, we have the
following classification result for complete gradient shrinking Ricci
solitons that have harmonic Weyl tensor. This extends the results from \cite%
{NW, CWZ, PW}.

\begin{thm}
\label{thm-class} Any $n-$dimensional complete gradient shrinking Ricci
soliton with harmonic Weyl tensor is a finite quotient of $\mathbb{R}^{n},$ $%
\mathbb{S}^{n-1}\times \mathbb{R}$ or $\mathbb{S}^{n}$.
\end{thm}

The fact that shrinking and steady solitons have many properties common to
manifolds with non-negative Ricci curvature motivates us to study the issue
of existence of harmonic functions on these manifolds. It is known, see \cite%
{LT}, that the existence of certain classes of harmonic functions is related
to the existence of ends of the manifold. Some other results about the
topology of shrinking Ricci solitons have been obtained in \cite{W, FMZ}.

We first recall some known terminology. A manifold is called nonparabolic if
it admits a positive symmetric Green's function. Otherwise it is called
parabolic. An end of a manifold is called nonparabolic if it admits a
positive symmetric Green's function that satisfies the Neumann boundary
condition on the boundary of the end. Otherwise, it is called parabolic.

We recall that a gradient shrinking K\"{a}hler-Ricci soliton satisfies 
\begin{equation*}
R_{\alpha \bar{\beta}}+f_{\alpha \bar{\beta}}=g_{\alpha \bar{\beta}},
\end{equation*}
for a smooth function $f$ that has the property $f_{\alpha \beta }=f_{\bar{%
\alpha}\bar{\beta}}=0.$ We will establish the following Liouville-type
theorem for gradient shrinking K\"{a}hler-Ricci solitons.

\begin{thm}
\label{thm-no-har}Let $(M,g)$ be a gradient shrinking K\"{a}hler-Ricci
soliton. If $u$ is a harmonic function with $\int_{M}|\nabla u|^{2}\,<\infty 
$ then $u$ is a constant function.
\end{thm}

Let us point out that on a manifold with a weighted volume $e^{-f}dv,$ a
naturally defined operator is the $f-$Laplacian $\Delta _{f}:=\Delta -\nabla
f\cdot \nabla ,$ which is self adjoint with respect to the weighted volume.
Though this operator will be used in our proofs, let us point out that in
Theorem \ref{thm-no-har} and everywhere in the paper the assumption of being
harmonic refers to the usual Laplace operator i.e. $\Delta u=0.$

As a consequence of the previous theorem we proved that any gradient
shrinking K\"{a}hler-Ricci soliton has at most one nonparabolic end.
Furthermore, if in addition to being K\"{a}hler we have an upper bound on
the scalar curvature of the form $\sup_{M}R<\frac{n}{2}-1$ on $M,$ then $M$
is connected at infinity, i.e., it has one end.

Let us now consider the case of gradient steady Ricci solitons $(M,g),$
which by definition satisfy the equation 
\begin{equation*}
Ric=Hess\left( f\right) .
\end{equation*}
We will prove the following results.

\begin{thm}
\label{thm-steady} Let $(M,g)$ be a gradient steady Ricci soliton. If $u$ is
harmonic with $\int_{M}|\nabla u|^{2}<\infty $ then $u$ is constant on $M$.
\end{thm}

As a consequence, we have that any gradient steady Ricci soliton has at most
one nonparabolic end. In the case we assume more on geometry of steady
solitons we can prove the following structural result.

\begin{thm}
\label{thm-steady-two-ends} If $M$ is a gradient steady K\"{a}hler-Ricci
soliton with Ricci curvature bounded below and such that for any $x\in M$, $%
Vol\left( B_{x}\left( 1\right) \right) \geq C>0$ for some constant $C$
independent of $x$, then either it is connected at infinity or it splits
isometrically as $M=\mathbb{R}\times N$, for a compact Ricci flat manifold $%
N $.
\end{thm}

In proving these results, a good knowledge of the volume growth and
asymptotics of the potential function is quite important. While for gradient
shrinking solitons the results  in \cite{CD} provide such estimates, not
much is known for steady Ricci solitons. In this sense, we have established
the following result.

\begin{thm}
\label{thm-inf-vol} Let $(M,g)$ be any gradient steady Ricci soliton. There
exist constants $a,c,r_{0}>0$ so that for any $r>r_{0}$, 
\begin{equation*}
ce^{a\sqrt{r}}\geq Vol(B_{p}(r))\geq c\cdot r.
\end{equation*}
\end{thm}

As a consequence of this Theorem, we can prove estimates for the potential
function $f,$ see Corollary \ref{potential}. Related to Theorem \ref%
{thm-inf-vol}, we note that it is known that shrinking Ricci solitons with
bounded Ricci curvature have at least linear volume growth, see \cite{CN}.
This is a consequence of the log Sobolev inequality established by Carillo
and Ni in \cite{CN} and Perelman's argument in \cite{P}.

The organization of the paper is as follows. In Section \ref{sec-integral}
we prove Theorem \ref{prop-int-ric0} and Theorem \ref{prop-iden0}. We show
how to use them to give the proof of Theorem \ref{thm-class}. In Section \ref%
{sec-kahler} we prove Theorem \ref{thm-no-har} about gradient shrinking K%
\"{a}hler-Ricci solitons. In Section \ref{sec-steady} we prove Theorem \ref%
{thm-steady} and Theorem \ref{thm-steady-two-ends}. In Section \ref%
{vol_steady} we prove Theorem \ref{thm-inf-vol} about volume of any gradient
steady Ricci soliton.

\section{Integral curvature estimates for shrinking solitons}

\label{sec-integral}

Let $(M,g)$ be a complete gradient shrinking Ricci soliton, given by
equation 
\begin{equation*}
R_{ij}+f_{ij}=\frac{1}{2}g_{ij}.
\end{equation*}%
It is well known that after normalizing the potential function we have the
following set of identities satisfied by the soliton, 
\begin{equation}
R+\Delta f=\frac{n}{2},\ \ \ \ \left\vert \nabla f\right\vert ^{2}+R=f\ \ \ 
\text{and \ \ \ }\nabla _{i}\left( R_{ij}e^{-f}\right) =0.
\label{properties}
\end{equation}%
We have denoted with $R$ the scalar curvature of $M.$ In \cite{CD} Cao and
Zhou have proved that for any fixed point $p\in M$ there are uniform
constants $C,c>0$ so that 
\begin{eqnarray}
Vol(B_{p}(r)) &\leq &Cr^{n}\,\,\,\mbox{and}  \label{eq-cao} \\
\frac{1}{4}\left( r\left( x\right) -c\right) ^{2} &\leq &f\left( x\right)
\leq \frac{1}{4}\left( r\left( x\right) +c\right) ^{2},  \notag
\end{eqnarray}%
where $r(x)=\mathrm{dist}(x,p)$. We note that asymptotic estimates for the
potential function $f$ were previously studied in \cite{P, Ni, NW, FMZ} and
interesting volume growth properties of solitons were also investigated in 
\cite{CN}.

In \cite{Ch} Chen proved that every complete ancient solution to the Ricci
flow has nonnegative scalar curvature for all times of their existence, see
also Proposition 5.5 in \cite{C1}. In particular, this holds for gradient
shrinking Ricci solitons and therefore using (\ref{eq-cao}) and (\ref%
{properties}) we have: 
\begin{equation}
0\leq R\leq \frac{(r(x)+c)^{2}}{4}.  \label{eq-scalar-bound}
\end{equation}

In this section we will prove Theorem \ref{prop-int-ric0} and Theorem \ref%
{prop-iden0}. We will say how to use them to prove Theorem \ref{thm-class}.

\begin{thm}
\label{prop-int-ric} For any complete gradient shrinking Ricci soliton $%
\left( M,g\right) $ we have 
\begin{equation*}
\int_{M}\left\vert Ric\right\vert ^{2}e^{-\lambda f}<\infty ,\ \ \text{for
any }\lambda >0.
\end{equation*}
\end{thm}

\begin{proof}[Proof of Theorem \protect\ref{prop-int-ric}]
For a cut-off function $\phi $ on $M$ we have, integrating by parts and
using (\ref{properties}): 
\begin{gather}
\int_{M}\left\vert Ric\right\vert ^{2}e^{-\lambda f}\phi
^{2}=\int_{M}R_{ij}\left( \frac{1}{2}g_{ij}-f_{ij}\right) e^{-\lambda f}\phi
^{2}  \label{p1} \\
=\frac{1}{2}\int_{M}Re^{-\lambda f}\phi ^{2}+\int_{M}f_{i}\nabla _{j}\left(
R_{ij}e^{-\lambda f}\phi ^{2}\right)  \notag \\
=\frac{1}{2}\int_{M}Re^{-\lambda f}\phi ^{2}+\left( 1-\lambda \right)
\int_{M}R_{ij}f_{i}f_{j}e^{-\lambda f}\phi
^{2}+\int_{M}R_{ij}f_{i}e^{-\lambda f}\left( \phi ^{2}\right) _{j}.  \notag
\end{gather}

By simple algebraic manipulations we have:%
\begin{eqnarray*}
\left( 1-\lambda \right) \int_{M}R_{ij}f_{i}f_{j}e^{-\lambda f}\phi ^{2}
&\leq &\frac{1}{4}\int_{M}\left\vert Ric\right\vert ^{2}e^{-\lambda f}\phi
^{2}+\left\vert 1-\lambda \right\vert ^{2}\int_{M}\left\vert \nabla
f\right\vert ^{4}e^{-\lambda f}\phi ^{2} \\
\int_{M}R_{ij}f_{i}e^{-\lambda f}\left( \phi ^{2}\right) _{j} &\leq &\frac{1%
}{4}\int_{M}\left\vert Ric\right\vert ^{2}e^{-\lambda f}\phi
^{2}+4\int_{M}\left\vert \nabla f\right\vert ^{2}e^{-\lambda f}\left\vert
\nabla \phi \right\vert ^{2}.
\end{eqnarray*}

Notice that from (\ref{eq-cao}) we know $\int_{M}\left\vert \nabla
f\right\vert ^{4}e^{-\lambda f}<\infty $ and $\int_{M}Re^{-\lambda f}<\infty
.$ Therefore, from (\ref{p1}) it easily follows that $\int_{M}\left\vert
Ric\right\vert ^{2}e^{-\lambda f}<\infty .$ This proves the Theorem.

Let us also point out that in the special case when $\lambda =1,$ (\ref{p1})
implies in particular that 
\begin{equation*}
\int_{M}\left\vert Ric\right\vert ^{2}e^{-f}=\frac{1}{2}\int_{M}Re^{-f}<%
\infty .
\end{equation*}
\end{proof}

Next denote 
\begin{eqnarray*}
\left\vert \nabla Ric\right\vert ^{2} &=&\sum \left\vert \nabla
_{k}R_{ij}\right\vert ^{2} \\
div\left( Rm\right) _{jkl} &=&\nabla _{i}R_{ijkl}
\end{eqnarray*}

\begin{thm}
\label{prop-iden} Let $(M,g)$ be a gradient shrinking Ricci soliton. If \
for some $\lambda <1$ we have $\int_{M}\left\vert Rm\right\vert
^{2}e^{-\lambda f}<\infty ,$ then the following identity holds 
\begin{equation*}
\int_{M}\left\vert \nabla Ric\right\vert ^{2}e^{-f}=\int_{M}\left\vert
div\left( Rm\right) \right\vert ^{2}e^{-f}<\infty .
\end{equation*}
\end{thm}

\begin{proof}[Proof of Theorem \protect\ref{prop-iden}]
We use the following formulas, true for any gradient shrinking Ricci
soliton: 
\begin{eqnarray}
\nabla _{k}R_{ij}-\nabla _{j}R_{ik} &=&R_{kjhi}f_{h}  \label{b} \\
\nabla _{i}\left( R_{ijkl}e^{-f}\right) &=&0  \label{a} \\
div\left( Rm\right) _{jkl} &=&R_{lkjp}f_{p}  \label{c}
\end{eqnarray}%
It is known that the Ricci curvature satisfies the equation 
\begin{equation*}
\Delta _{f}R_{ij}=R_{ij}-2R_{ikjl}R_{kl}
\end{equation*}%
were $\Delta _{f}R_{ij}:=\Delta R_{ij}-\langle \nabla f,\nabla R_{ij}\rangle 
$ is the $f-$Laplacian of the Ricci tensor.

For a cut-off function $\phi $ on $M$ we have 
\begin{gather*}
\int_{M}\left\vert \nabla Ric\right\vert ^{2}e^{-f}\phi ^{2}=-\int_{M}\left(
\Delta _{f}R_{ij}\right) R_{ij}e^{-f}\phi ^{2}-\int_{M}\left( \nabla
_{k}R_{ij}\right) R_{ij}e^{-f}\left( \phi ^{2}\right) _{k} \\
=-\int_{M}\left\vert R_{ij}\right\vert ^{2}e^{-f}\phi
^{2}+2\int_{M}R_{ikjl}R_{ij}R_{kl}e^{-f}\phi ^{2}-\int_{M}\left( \nabla
_{k}R_{ij}\right) R_{ij}e^{-f}\left( \phi ^{2}\right) _{k}
\end{gather*}%
The Riemann curvature term can be computed using the soliton equation: 
\begin{equation*}
2\int_{M}R_{ikjl}R_{ij}R_{kl}e^{-f}\phi ^{2}=\int_{M}\left\vert
R_{ij}\right\vert ^{2}e^{-f}\phi
^{2}-2\int_{M}R_{ikjl}R_{ij}f_{kl}e^{-f}\phi ^{2}.
\end{equation*}%
This gives 
\begin{equation*}
\int_{M}\left\vert \nabla Ric\right\vert ^{2}e^{-f}\phi
^{2}=-2\int_{M}R_{ikjl}R_{ij}f_{kl}e^{-f}\phi ^{2}-\int_{M}\left( \nabla
_{k}R_{ij}\right) R_{ij}e^{-f}\left( \phi ^{2}\right) _{k}.
\end{equation*}%
Using (\ref{a}) we now get 
\begin{gather*}
-2\int_{M}R_{ikjl}R_{ij}f_{kl}e^{-f}\phi ^{2}=2\int_{M}f_{k}\nabla
_{l}\left( R_{ikjl}e^{-f}R_{ij}\phi ^{2}\right) \\
=2\int_{M}f_{k}R_{ikjl}e^{-f}\nabla _{l}\left( R_{ij}\phi ^{2}\right) \\
=2\int_{M}R_{ikjl}\left( \nabla _{l}R_{ij}\right) f_{k}e^{-f}\phi
^{2}+2\int_{M}R_{ikjl}R_{ij}f_{k}e^{-f}\left( \phi ^{2}\right) _{l}.
\end{gather*}%
Notice moreover, using (\ref{b}) and (\ref{c}) that 
\begin{gather*}
2R_{ikjl}\left( \nabla _{l}R_{ij}\right) f_{k}=-2R_{ljik}f_{k}\left( \nabla
_{l}R_{ij}\right) =2R_{ljik}f_{k}\left( \nabla _{j}R_{il}\right) \\
=\left( R_{ljik}f_{k}\right) \left( \nabla _{j}R_{il}-\nabla
_{l}R_{ij}\right) =\left\vert divRm\right\vert ^{2}
\end{gather*}%
This proves that 
\begin{gather}
\int_{M}\left\vert \nabla Ric\right\vert ^{2}e^{-f}\phi
^{2}=\int_{M}\left\vert divRm\right\vert ^{2}e^{-f}\phi ^{2}  \label{5} \\
+2\int_{M}R_{ikjl}R_{ij}f_{k}e^{-f}\left( \phi ^{2}\right)
_{l}-\int_{M}\left( \nabla _{k}R_{ij}\right) R_{ij}e^{-f}\left( \phi
^{2}\right) _{k}.  \notag
\end{gather}%
Our hypothesis and soliton identities imply that 
\begin{eqnarray*}
\int_{M}\left\vert divRm\right\vert ^{2}e^{-f} &\leq &C\int |Rm|^{2}|\nabla
f|^{2}e^{-f}\leq C\int_{M}|Rm|^{2}e^{-\lambda f}<\infty \\
\int_{M}\left\vert R_{ikjl}R_{ij}f_{k}\phi _{l}\right\vert e^{-f} &\leq
&C\int_{M}|Rm|^{2}|\nabla f|e^{-f}\leq C\int_{M}|Rm|^{2}e^{-\lambda
f}<\infty ,
\end{eqnarray*}%
for $\lambda <1$ given in the statement of Theorem \ref{prop-iden}. Hence, (%
\ref{5}) and the arithmetic-mean inequality imply that 
\begin{eqnarray*}
\int_{M}\left\vert \nabla Ric\right\vert ^{2}e^{-f}\phi ^{2} &\leq
&C+2\int_{M}\left\vert \nabla _{k}R_{ij}\right\vert \left\vert
R_{ij}\right\vert e^{-f}\phi \left\vert \nabla \phi \right\vert \\
&\leq &C+\frac{1}{2}\int_{M}\left\vert \nabla Ric\right\vert ^{2}e^{-f}\phi
^{2}+2\int_{M}\left\vert Ric\right\vert ^{2}e^{-f}\left\vert \nabla \phi
\right\vert ^{2}.
\end{eqnarray*}%
This clearly shows that 
\begin{equation*}
\int_{M}\left\vert \nabla Ric\right\vert ^{2}e^{-f}<\infty .
\end{equation*}
Returning to (\ref{5}) we see that all terms involving $\nabla \phi $ must
converge to zero as $r\rightarrow \infty .$ More precisely, taking $\phi $
such that $\phi =1$ on $B_{p}\left( r\right) ,$ $\phi =0$ on $M\backslash
B_{p}\left( 2r\right) $ and $\left\vert \nabla \phi \right\vert \leq \frac{c%
}{r}$, it follows that as $r\rightarrow \infty $ 
\begin{eqnarray*}
\left\vert \int_{M}R_{ikjl}R_{ij}f_{k}e^{-f}(\phi ^{2})_{l}\right\vert &\leq
&C\int_{B_{p}(2r)\backslash B_{p}(r)}|Rm|^{2}e^{-\lambda f}\rightarrow 0, \\
\left\vert \int_{M}(\nabla _{k}R_{ij})R_{ij}e^{-f}(\phi ^{2})_{k}\right\vert
&\leq &C\left( \int_{M}|\nabla Ric|^{2}e^{-f}\right) ^{\frac{1}{2}}\left(
\int_{B_{p}(2r)\backslash B_{p}(r)}|Ric|^{2}e^{-f}\right) ^{\frac{1}{2}%
}\rightarrow 0.
\end{eqnarray*}%
This proves the Theorem.
\end{proof}

\begin{remark}
In \cite{CWZ} the integral identity in Theorem \ref{prop-iden} was
established under a pointwise assumption on the Riemann tensor, that is, 
\begin{equation*}
\left\vert R_{ijkl}\right\vert \left( x\right) \leq e^{ar\left( x\right) +1}.
\end{equation*}%
Then Shi's derivative estimate will imply that $\left\vert \nabla
Ric\right\vert $ has a similar growth and therefore $\int_{M}\left\vert
\nabla Ric\right\vert ^{2}e^{-f}<\infty .$ Then it follows that the
integration by parts argument is valid in the noncompact setting. The
advantage of our argument is that it requires only weak integral control of
the Riemann tensor.
\end{remark}

\begin{cor}
\label{cor-iden} Let $\left( M,g\right) $ be a gradient shrinking Ricci
soliton with harmonic Weyl tensor. Then%
\begin{equation*}
\int_{M}\left\vert \nabla Ric\right\vert ^{2}e^{-f}=\int_{M}\left\vert div
\left( Rm\right) \right\vert ^{2}e^{-f}<\infty .
\end{equation*}
\end{cor}

\begin{proof}[Proof of Corollary \protect\ref{cor-iden}]
We start with a formula established in in the proof of Theorem \ref%
{prop-iden}: 
\begin{gather*}
\int_{M}\left\vert \nabla Ric\right\vert ^{2}e^{-f}\phi
^{2}=\int_{M}\left\vert divRm\right\vert ^{2}e^{-f}\phi ^{2} \\
+2\int_{M}R_{ikjl}R_{ij}f_{k}e^{-f}\left( \phi ^{2}\right)
_{l}-\int_{M}\left( \nabla _{k}R_{ij}\right) R_{ij}e^{-f}\left( \phi
^{2}\right) _{k},
\end{gather*}%
It is known that if the Weyl tensor is harmonic i.e. $divW=0$, we have the
following identity for gradient shrinkers, \cite{FG}: 
\begin{equation}
R_{ijkl}\nabla _{l}f=\frac{1}{\left( n-1\right) }%
(R_{il}f_{l}g_{jk}-R_{jl}f_{l}g_{ik}).  \label{c1}
\end{equation}%
Since on shrinking soliton we have the identity $%
(divRm)_{kji}=R_{ijkl}f_{l}, $ by (\ref{c1}) we obtain 
\begin{equation*}
\int_{M}|divRm|^{2}e^{-f}\leq C\int_{M}|Ric|^{2}|\nabla f|^{2}e^{-f}\leq
\int_{M}|Ric|^{2}e^{-\mu f}<\infty ,
\end{equation*}%
for $\mu <1$. Moreover, for a cut-off as in Theorem \ref{prop-iden} we get: 
\begin{eqnarray*}
\int_{M}|R_{ikjl}f_{k}R_{ij}(\phi ^{2})_{l}|e^{-f} &\leq &\frac{c}{r}\left(
\int_{M}|divRm|^{2}e^{-f}+\int_{M}|Ric|^{2}e^{-f}\right) \\
&\leq &\frac{C}{r}\,\,\rightarrow 0\,\,\,\mbox{as}\,\,\,r\rightarrow \infty .
\end{eqnarray*}%
The rest of the proof is same as the proof of Theorem \ref{prop-iden}.
\end{proof}

Either Theorem \ref{prop-iden} combined with results in \cite{CWZ} or
Theorem \ref{prop-int-ric} combined with the results in \cite{PW} can be
used to show that any locally conformally flat gradient shrinking Ricci
soliton is rigid. More generally, in the case of a harmonic Weyl tensor we
have the following.

\begin{thm}
\label{weyl}Any $n-$dimensional complete shrinking gradient Ricci soliton
with harmonic Weyl tensor is a finite quotient of $\mathbb{R}^{n},$ $\mathbb{%
S}^{n-1}\times \mathbb{R}$ or $\mathbb{S}^{n}$.
\end{thm}

\begin{proof}[Proof of Theorem \protect\ref{weyl}]
In \cite{FG} it has been proved that once we have Corollary \ref{cor-iden},
then any gradient shrinking Ricci soliton with harmonic Weyl tensor must be
a finite quotient as in the statement of the Theorem. This proves the
Theorem.
\end{proof}

\section{Gradient shrinking Ricci solitons and harmonic functions}

\label{sec-kahler}

As mentioned in the introduction, shrinking solitons have many properties in
common with manifolds with non-negative Ricci curvature. For manifolds with
non-negative Ricci curvature S.-T. Yau (\cite{Y}) proved that positive
harmonic functions are necessarily trivial. The question is whether this
generalizes to gradient shrinking Ricci solitons. One motivation for
studying the existence of harmonic functions comes from its relation to the
structure of the manifold at infinity, that is, the number of ends.

Let $(M,g)$ be a gradient shrinking K\"{a}hler-Ricci soliton, that is, 
\begin{equation*}
R_{\alpha \bar{\beta}}+f_{\alpha \bar{\beta}}=g_{\alpha \bar{\beta}},
\end{equation*}
for a smooth function $f$ that has the property 
\begin{equation}
f_{\alpha \beta }=f_{\bar{\alpha}\bar{\beta}}=0.  \label{eq-kr-pot}
\end{equation}%
Note that in complex coordinates we have, for any $u$, $v\in C^{\infty
}\left( M\right) $ 
\begin{gather*}
\langle \nabla u,\nabla v\rangle =\frac{1}{2}\left( u_{\alpha }v_{\bar{\alpha%
}}+u_{\bar{\alpha}}v_{\alpha }\right) \\
\Delta u=u_{\alpha \bar{\alpha}}
\end{gather*}%
Our next result says there are no harmonic functions with bounded total
energy on complete gradient shrinking K\"{a}hler-Ricci solitons.

\begin{thm}
\label{thm-no-har0} Let $(M,g)$ be a gradient shrinking K\"{a}hler-Ricci
soliton. If $u$ is a harmonic function with $\int_{M}|\nabla u|^{2}\,<\infty 
$ then $u$ is a constant function.
\end{thm}

\begin{proof}[Proof of Theorem \protect\ref{thm-no-har0}]
Let $u$ satisfy $\Delta u=0\ $on $M$ and\ $\int_{M}|\nabla u|^{2}\,<\infty .$
We first prove that $\nabla f$ and $\nabla u$ are orthogonal to each other.
Then this implies that $u$ is in fact $f-$harmonic and this fact forces $u$
to be constant.

Let $\phi :M\rightarrow \lbrack 0,1]$ be a cut off function such that $\phi
=1$ on $B_{p}\left( r\right) $ (a geodesic ball centered at some fixed point 
$p\in M$ of radius $r$), $\phi =0$ outside $B_{p}\left( 2r\right) $ and $%
|\nabla \phi |\leq \frac{C}{r}$.

We recall that, according to a result of P. Li in \cite{L}, if $u$ is
harmonic and with finite total energy on a K\"{a}hler manifold then it is in
fact pluriharmonic, that is $u_{\alpha \bar{\beta}}=0.$

Let us define $F\in C^{\infty }(M)$ to be 
\begin{equation*}
F:=\langle \nabla f,\nabla u\rangle =\frac{1}{2}(u_{\alpha }f_{\bar{\alpha}%
}+u_{\bar{\alpha}}f_{\alpha }).
\end{equation*}%
We show that $F\equiv 0$. To this end, observe that 
\begin{equation*}
(u_{\alpha }f_{\bar{\alpha}})_{\bar{\delta}}=u_{\alpha \bar{\delta}}f_{\bar{%
\alpha}}+u_{\alpha }f_{\bar{\alpha}\bar{\delta}}=0,
\end{equation*}%
where we have used (\ref{eq-kr-pot}) and $u_{\alpha \bar{\beta}}=0$.
Similarly, $(u_{\bar{\alpha}}f_{\alpha })_{\delta }=0$. This implies that 
\begin{equation*}
\Delta F=0,
\end{equation*}%
hence we have:%
\begin{gather*}
\int_{M}|\nabla F|^{2}\phi ^{2}\,=-\int_{M}\left( \Delta F\right) {F}\phi
^{2}\,-2\int_{M}{F}\phi \langle \nabla F,\nabla \phi \rangle  \\
\leq 2\int_{M}|\nabla F||F||\nabla \phi |\phi \,\leq \frac{1}{2}%
\int_{M}|\nabla F|^{2}\phi ^{2}+2\int_{M}|F|^{2}|\nabla \phi |^{2}.
\end{gather*}%
From here we get 
\begin{eqnarray*}
\int_{M}|\nabla F|^{2}\phi ^{2}\, &\leq &4\int_{M}|F|^{2}|\nabla \phi
|^{2}\leq 4\int_{M}|\nabla u|^{2}|\nabla f|^{2}|\nabla \phi |^{2} \\
&\leq &C\int_{B_{p}\left( 2r\right) \backslash B_{p}\left( r\right) }|\nabla
u|^{2}\rightarrow 0\,\,\,\mbox{as}\,\,\,r\rightarrow \infty .
\end{eqnarray*}%
In the last line we have used that $\left\vert \nabla \phi \right\vert \leq 
\frac{C}{r},$ the fact that $\left\vert \nabla f\right\vert $ grows linearly
on $M$ and the assumption that $u$ has finite total energy. This yields $F=%
\mathrm{const}$ on $M$. The asymptotic behavior (\ref{eq-cao}) of $f$
guarantees that $f$ attains its minimum somewhere on a compact subset of $M$
and therefore its gradient vanishes at the minimum point. In particular,
that means $F\equiv 0$ on $M$. Define the $f$-Laplacian of a function to be 
\begin{equation*}
\Delta _{f}=\Delta -\nabla f\cdot \nabla .
\end{equation*}%
We now prove that we have the following:%
\begin{equation}
\Delta _{f}|\nabla u|\geq \frac{1}{2}\left\vert \nabla u\right\vert ,\ \ 
\text{whenever }\left\vert \nabla u\right\vert \neq 0.  \label{f-laplace}
\end{equation}%
Since we have proved $\langle \nabla u,\nabla f\rangle =0$, it follows that 
\begin{equation*}
\Delta _{f}u=\Delta u-\left\langle \nabla f,\nabla u\right\rangle =0.
\end{equation*}%
The Bochner formula implies%
\begin{eqnarray*}
\Delta _{f}\left\vert \nabla u\right\vert ^{2} &=&2Ric_{f}\left( \nabla
u,\nabla u\right) +2\left\langle \nabla \Delta _{f}u,\nabla u\right\rangle
+2\left\vert u_{ij}\right\vert ^{2} \\
&=&\left\vert \nabla u\right\vert ^{2}+2\left\vert u_{ij}\right\vert
^{2}\geq \left\vert \nabla u\right\vert ^{2}+2\left\vert \nabla \left\vert
\nabla u\right\vert \right\vert ^{2}.
\end{eqnarray*}%
In the last line we have used the Kato inequality. Since on the other hand, 
\begin{equation*}
\Delta _{f}\left\vert \nabla u\right\vert ^{2}=2\left\vert \nabla
u\right\vert \Delta _{f}\left\vert \nabla u\right\vert +2\left\vert \nabla
\left\vert \nabla u\right\vert \right\vert ^{2},
\end{equation*}%
it is clear that we get (\ref{f-laplace}).

Since in particular, $\Delta _{f}|\nabla u|\geq 0$ and $\int_{M}\left\vert
\nabla u\right\vert ^{2}e^{-f}<\infty $, from a standard argument of Yau, 
see also \cite{N} and Theorem 4.2 in \cite{PW} for the $f-$Laplacian case,
it follows that $| \nabla u|=C$ on $M$. Then (\ref{f-laplace}) implies that $%
| \nabla u|=0$, hence $u$ is constant on $M$.
\end{proof}

\begin{cor}
\label{cor-ends} Let $M$ be a gradient shrinking K\"{a}hler-Ricci soliton.
It then has at most one nonparabolic end.
\end{cor}

\begin{proof}[Proof of Corollary \protect\ref{cor-ends}]
From the theory of Li and Tam \cite{LT} it is known that if a manifold $M$
has at least two nonparabolic ends, then there exists a bounded harmonic
function $u$ on $M$ which has finite total energy $\int_{M}$ $\left\vert
\nabla u\right\vert ^{2}<\infty .$ This is impossible by Theorem \ref%
{thm-no-har0}.
\end{proof}

The next result shows that under some upper bound for the scalar curvature
of $M,$ all ends are nonparabolic. Here we do not make the assumption of $M$
being K\"{a}hler.

\begin{prop}
\label{prop-general} Let $(M,g)$ be a gradient shrinking Ricci soliton such
that for some constant $\alpha $ we have $R\leq \alpha <\frac{n}{2}-1$. Then
all the ends of $M$ are nonparabolic.
\end{prop}

\begin{proof}[Proof of Proposition \protect\ref{prop-general}]
For $a=\frac{n}{2}-\alpha -1>0$ we compute 
\begin{eqnarray*}
\Delta f^{-a} &=&-af^{-a-1}\Delta f+a(a+1)\left\vert \nabla f\right\vert
^{2}f^{-a-2} \\
&=&\left( -a\left( \frac{n}{2}-R\right) +a(a+1)\right) f^{-a-1}-a\left(
a+1\right) Rf^{-a-2} \\
&\leq &a\left( \alpha -\frac{n}{2}+a+1\right) f^{-a-1}=0.
\end{eqnarray*}%
This proves that there exists a positive super-harmonic function which
converges to zero at infinity. Then it is known \cite{L1, G} that any end of 
$M$ (and hence $M$) is nonparabolic.
\end{proof}

\begin{remark}
If $R=\frac{n}{2}-1$ the conclusion in Proposition \ref{prop-general} no
longer holds. For example, $\mathbb{S}^{n-2}\times \mathbb{R}^{2},$ where $%
\mathbb{R}^{2}$ is the Gaussian soliton has $R=\frac{n}{2}-1$ and it is
parabolic.
\end{remark}

Combining Proposition \ref{prop-general} and Corollary \ref{cor-ends} we
conclude the following.

\begin{cor}
\label{nonp-ends} If $(M,g)$ is a gradient shrinking K\"{a}hler-Ricci
soliton with scalar curvature 
\begin{equation*}
R\leq \alpha <\frac{n}{2}-1
\end{equation*}%
for some constant $\alpha $, then it is connected at infinity.
\end{cor}

We conclude this section with the observation that under the hypothesis in
Proposition \ref{prop-general} the manifold satisfies a Poincar\'{e}
inequality, similar to Hardy's inequality. Indeed, for any function $\phi $
with compact support, 
\begin{gather*}
\int_{M}\left( \Delta f\right) f^{-1}\phi ^{2}=\int_{M}f^{-2}\left\vert
\nabla f\right\vert ^{2}\phi ^{2}-2\int_{M}\phi f^{-1}\langle \nabla
f,\nabla \phi \rangle \\
\leq \left( 1+\varepsilon \right) \int_{M}f^{-2}\left\vert \nabla
f\right\vert ^{2}\phi ^{2}+\varepsilon ^{-1}\int_{M}\left\vert \nabla \phi
\right\vert ^{2}.
\end{gather*}

On the other hand, 
\begin{equation*}
\int_{M}\left( \Delta f\right) f^{-1}\phi ^{2}\geq \left( \frac{n}{2}-\alpha
\right) \int_{M}f^{-1}\phi ^{2}\ \ \text{and\ }\int_{M}f^{-2}\left\vert
\nabla f\right\vert ^{2}\phi ^{2}\leq \int_{M}f^{-1}\phi ^{2}.
\end{equation*}%
Therefore, we arrive at 
\begin{equation*}
\varepsilon \left( \frac{n}{2}-\alpha -1-\varepsilon \right)
\int_{M}f^{-1}\phi ^{2}\leq \int_{M}\left\vert \nabla \phi \right\vert ^{2}.
\end{equation*}%
Choosing $\varepsilon =\frac{1}{2}\left( \frac{n}{2}-\alpha -1\right) >0$ we
obtain a weighted Poincar\'{e} inequality 
\begin{equation*}
\int_{M}\rho \phi ^{2}\leq \int_{M}\left\vert \nabla \phi \right\vert ^{2},\
\ \text{for \ }\rho :=\frac{1}{4}\left( \frac{n}{2}-\alpha -1\right)
^{2}f^{-1}.
\end{equation*}

Weighted Poincar\'{e} inequalities are known to be equivalent to the
manifold being nonparabolic \cite{LW}. Moreover, $M$ satisfies the property $%
\left( P_\rho \right) $ in the sense of \cite{LW}. For gradient steady
solitons a weighted Poincar\'{e} inequality was proved in \cite{CN}.

\section{Steady gradient Ricci solitons}

\label{sec-steady}

In this section we will study the existence of harmonic functions on
gradient steady solitons, which is, as we have mentioned earlier, tightly
related to the structure of a given manifold.

Let $(M,g)$ be a gradient steady Ricci soliton, that is, 
\begin{equation*}
Ric=Hess\left( f\right) .
\end{equation*}
In \cite{Ha2} Hamilton showed that 
\begin{equation*}
R+|\nabla f|^{2}=\lambda ,
\end{equation*}%
for some constant $\lambda >0$. Every steady soliton is in particular an
ancient solution to the Ricci flow and therefore $R\geq 0$. This implies $%
|\nabla f|\leq \sqrt{\lambda }$ and therefore 
\begin{equation}
f(x)\leq f(p)+\sqrt{\lambda }r(x),  \label{eq-lin}
\end{equation}%
where $p\in M$ is a fixed point, $r(x)=\mathrm{dist}(x,p)$ and $x\in M$ is
an arbitrary point on $M$. We claim the following result.

\begin{thm}
\label{th-steady} If $(M,g)$ is a gradient steady Ricci soliton and $\Delta
u=0$ with $\int_{M}|\nabla u|^{2}<\infty $, then $u$ is a constant function.
\end{thm}

\begin{proof}[Proof of Theorem \protect\ref{th-steady}]
For a cut-off $\phi $ on $M$ it follows that 
\begin{equation}
\int_{M}Ric\left( \nabla u,\nabla u\right) \phi
^{2}=\int_{M}f_{ij}u_{i}u_{j}\phi ^{2}=-\int_{M}u_{ij}f_{i}u_{j}\phi
^{2}-\int_{M}f_{i}u_{i}u_{j}\left( \phi ^{2}\right) _{j}.  \label{d}
\end{equation}%
Note that, integrating by parts we have:%
\begin{equation*}
-\int_{M}u_{ij}f_{i}u_{j}\phi ^{2}=\frac{1}{2}\int_{M}\left( \Delta f\right)
\left\vert \nabla u\right\vert ^{2}\phi ^{2}+\frac{1}{2}\int_{M}\left\vert
\nabla u\right\vert ^{2}\left\langle \nabla f,\nabla \phi ^{2}\right\rangle .
\end{equation*}%
Plug this in formula (\ref{d}) and get that 
\begin{eqnarray}
\int_{M}Ric\left( \nabla u,\nabla u\right) \phi ^{2} &=&\frac{1}{2}%
\int_{M}R\left\vert \nabla u\right\vert ^{2}\phi ^{2}+\frac{1}{2}%
\int_{M}\left\vert \nabla u\right\vert ^{2}\langle \nabla f,\nabla \phi
^{2}\rangle  \label{Ric} \\
&-&\int_{M}\langle \nabla f,\nabla u\rangle \cdot \langle \nabla u,\nabla
\phi ^{2}\rangle .  \notag
\end{eqnarray}%
A similar integration by parts argument was used in \cite{DX}. We now recall
the Bochner formula 
\begin{equation*}
\Delta \left\vert \nabla u\right\vert ^{2}=2Ric\left( \nabla u,\nabla
u\right) +2\left\vert u_{ij}\right\vert ^{2}\geq 2Ric\left( \nabla u,\nabla
u\right) +2\left\vert \nabla \left\vert \nabla u\right\vert \right\vert ^{2}.
\end{equation*}%
We multiply this by $\phi ^{2},$ use (\ref{Ric}), and integrate by parts: 
\begin{gather*}
2\int_{M}\left\vert \nabla \left\vert \nabla u\right\vert \right\vert
^{2}\phi ^{2}+\int_{M}R\left\vert \nabla u\right\vert ^{2}\phi ^{2}\leq
-\int_{M}\langle \nabla \left\vert \nabla u\right\vert ^{2},\nabla \phi
^{2}\rangle \\
-\int_{M}\left\vert \nabla u\right\vert ^{2}\langle \nabla f,\nabla \phi
^{2}\rangle +2\int_{M}\langle \nabla f,\nabla u\rangle \cdot \langle \nabla
u,\nabla \phi ^{2}\rangle \\
\leq \int_{M}\left\vert \nabla \left\vert \nabla u\right\vert \right\vert
^{2}\phi ^{2}+4\int_{M}\left\vert \nabla u\right\vert ^{2}\left\vert \nabla
\phi \right\vert ^{2} \\
-\int_{M}\left\vert \nabla u\right\vert ^{2}\langle \nabla f,\nabla \phi
^{2}\rangle +2\int_{M}\langle \nabla f,\nabla u\rangle \cdot \langle \nabla
u,\nabla \phi ^{2}\rangle .
\end{gather*}%
We have thus proved that 
\begin{gather}
\int_{M}\left\vert \nabla \left\vert \nabla u\right\vert \right\vert
^{2}\phi ^{2}+\int_{M}R\left\vert \nabla u\right\vert ^{2}\phi ^{2}\leq
4\int_{M}\left\vert \nabla u\right\vert ^{2}\left\vert \nabla \phi
\right\vert ^{2}  \label{inequality} \\
-\int_{M}\left\vert \nabla u\right\vert ^{2}\langle \nabla f,\nabla \phi
^{2}\rangle +2\int_{M}\langle \nabla f,\nabla u\rangle \cdot \langle \nabla
u,\nabla \phi ^{2}\rangle  \notag \\
\leq C\int_{M}\left\vert \nabla u\right\vert ^{2}\left\vert \nabla \phi
\right\vert ,  \notag
\end{gather}%
where in the last line we have used that $\left\vert \nabla f\right\vert
\leq C$. Letting $r\rightarrow \infty $ and using that $u$ has finite total
energy it results that $\left\vert \nabla \left\vert \nabla u\right\vert
\right\vert =R\left\vert \nabla u\right\vert ^{2}=0.$ This implies $|\nabla
u|=C$. But since $M$ is nonparabolic, we know it has infinite volume. This
is in fact true in general, see Theorem \ref{thm-inf-vol}. But $%
\int_{M}\left\vert \nabla u\right\vert ^{2}<\infty ,$ therefore $\left\vert
\nabla u\right\vert =0$ and this proves the Theorem.
\end{proof}

We have the analogous result to Corollary \ref{cor-ends} in the case of
gradient steady Ricci solitons.

\begin{cor}
\label{cor-steady} Let $M$ be a gradient steady Ricci soliton. Then it has
at most one nonparabolic end.
\end{cor}

\begin{proof}[Proof of Corollary \protect\ref{cor-steady}]
As in the proof of Corollary \ref{cor-ends}, we apply the results in \cite%
{LT}. It is known that if a manifold $M$ has at least two nonparabolic ends,
then there exists a bounded harmonic function $u$ on $M$ which has finite
total energy $\int_{M}$ $\left\vert \nabla u\right\vert ^{2}<\infty .$ This
is impossible by Theorem \ref{th-steady}.
\end{proof}

If we assume more on geometry of $(M,g)$ we can say more about its ends.
Notice that the previous Corollary does not tell us anything about parabolic
ends if any. We will prove Theorem \ref{thm-steady-two-ends} in two steps,
depending whether $M$ is nonparabolic or parabolic, in the following two
propositions.

\begin{prop}
\label{prop-first-step} Assume $(M,g)$ is a complete, nonparabolic, gradient
steady K\"{a}hler-Ricci soliton with Ricci curvature bounded below and such
that for every $x\in M$, $Vol(B_{x}(1))\geq C$, for a uniform constant $C>0.$
Then $M$ is connected at infinity.
\end{prop}

\begin{proof}[Proof of Proposition \protect\ref{prop-first-step}]
Since the manifold is nonparabolic, it has at least one nonparabolic end.
Assume that $M$ has more than one end. By Corollary \ref{cor-steady} we may
assume it has a parabolic end, call it $H$. Then $E:=M\backslash H$ is
nonparabolic (since otherwise our manifold would be parabolic). By \cite{Na}
(see also \cite{NR}) there exists a positive harmonic function $u$ on $M$ so
that

\begin{enumerate}
\item[(i)] $\int_{E}|\nabla u|^{2}<\infty $, \ \ $inf_{E}u=0$ and

\item[(ii)] $\lim_{x\rightarrow \infty ,x\in H}u(x)=\infty $.
\end{enumerate}

We will obtain a contradiction by following a similar argument as in Theorem %
\ref{thm-no-har0}. We first show that $u$ is pluriharmonic and use it to
deduce that $\langle \nabla u,\nabla f\rangle =0$.

Let us start with the observation that there exists a uniform constant $C$
so that 
\begin{equation*}
\sup_{H}|\nabla u|\leq C.
\end{equation*}%
Indeed, this was proved in \cite{NR}, Theorem 2.1. Now we prove that there
exists a uniform constant $C$ so that 
\begin{equation}
\int_{B_{p}(r)}|\nabla u|^{2}\leq Cr.  \label{i1}
\end{equation}%
Since $\left\vert \nabla u\right\vert $ is bounded on $M,$ we have that if $%
x\in H$ and $r(x)\leq r$ then $u(x)\leq C\cdot r$. Then, using the co-area
formula it follows that 
\begin{eqnarray*}
\int_{B_{p}(r)}|\nabla u|^{2} &=&\int_{B_{p}(r)\cap E}|\nabla
u|^{2}+\int_{B_{p}(r)\cap H}|\nabla u|^{2}\leq C+\int_{\{u\leq C\cdot
r\}\cap H}|\nabla u|^{2} \\
&=&C+\int_{0}^{Cr}\int_{\{u=t\}\cap H}|\nabla u|\leq Cr.
\end{eqnarray*}%
Note that since $u$ is harmonic it follows $\int_{u=t}|\nabla u|=\mathrm{%
const}$.

Lemma 3.1 in \cite{L} and (\ref{i1}) now imply that $u$ is pluriharmonic. As
in Theorem \ref{thm-no-har0}, let us denote 
\begin{equation*}
F:=\left\langle \nabla f,\nabla u\right\rangle .
\end{equation*}%
Then, following the argument in Theorem \ref{thm-no-har0}, we get 
\begin{equation*}
\int_{M}|\nabla F|^{2}\phi ^{2}\leq \frac{C}{r^{2}}\cdot
\int_{B_{p}(2r)\backslash B_{p}(r)}|\nabla u|^{2}|\nabla f|^{2}\leq \frac{C}{%
r}\rightarrow 0\,\,\,\mbox{as}\,\,\,r\rightarrow \infty .
\end{equation*}%
We have used that $|\nabla f|\leq C$ in the case of a steady soliton, and (%
\ref{i1}). This implies that $F$ is constant on $M$ i.e., $\left\langle
\nabla u,\nabla f\right\rangle =a\in \mathbb{R}$. Moreover, we can show that
in fact $a=0$, because otherwise, $|a|=|\langle \nabla u,\nabla f\rangle
|\leq |\nabla u|\cdot |\nabla f|\leq C|\nabla u|$, which implies $|\nabla
u|\geq \delta >0$ on $M$. Since $\int_{E}|\nabla u|^{2}<\infty $ and any
nonparabolic end $E$ has infinite volume \cite{L1} , we get a contradiction.
This proved indeed that 
\begin{equation*}
\left\langle \nabla u,\nabla f\right\rangle \equiv 0.
\end{equation*}%
Recall the inequality (\ref{inequality}) obtained in the proof of Theorem %
\ref{th-steady}. 
\begin{gather}
\int_{M}\left\vert \nabla \left\vert \nabla u\right\vert \right\vert
^{2}\phi ^{2}+\int_{M}R\left\vert \nabla u\right\vert ^{2}\phi ^{2}\leq
4\int_{M}\left\vert \nabla u\right\vert ^{2}\left\vert \nabla \phi
\right\vert ^{2}  \label{inequality1} \\
-\int_{M}\left\vert \nabla u\right\vert ^{2}\langle \nabla f,\nabla \phi
^{2}\rangle +2\int_{M}\langle \nabla f,\nabla u\rangle \cdot \langle \nabla
u,\nabla \phi ^{2}\rangle .  \notag
\end{gather}%
This holds true here as well, since $u$ is harmonic. Let us choose the
cut-off $\phi $ as follows. On the nonparabolic end $E$ we define for $A$
large enough 
\begin{equation*}
\phi =\left\{ 
\begin{array}{c}
1 \\ 
A+1-r \\ 
0%
\end{array}%
\right. 
\begin{array}{l}
\text{on} \\ 
\text{on} \\ 
\text{on}%
\end{array}%
\begin{array}{l}
B_{p}\left( A\right) \cap E, \\ 
\left( B_{p}\left( A+1\right) \backslash B_{p}\left( A\right) \right) \cap E,
\\ 
E\backslash B_{p}\left( A+1\right) .%
\end{array}%
\end{equation*}%
On the parabolic end $H$ we define for $T$ large enough 
\begin{equation*}
\phi =\left\{ 
\begin{array}{c}
1 \\ 
\frac{2T-u}{T} \\ 
0%
\end{array}%
\right. 
\begin{array}{l}
\text{on} \\ 
\text{on} \\ 
\text{on}%
\end{array}%
\begin{array}{l}
u\leq T, \\ 
T<u<2T, \\ 
2T\leq u.%
\end{array}%
\end{equation*}%
Notice that $\phi $ defined in this manner is indeed with compact support,
as $u$ is proper on the parabolic end. Observe that since $\langle \nabla
u,\nabla f\rangle =0$ we have 
\begin{eqnarray*}
\int_{M}\langle \nabla f,\nabla u\rangle \cdot \langle \nabla u,\nabla \phi
^{2}\rangle  &=&0, \\
\int_{H}\left\vert \nabla u\right\vert ^{2}\langle \nabla f,\nabla \phi
^{2}\rangle  &=&0.
\end{eqnarray*}%
Notice that this is true for any $A$ and $T.$ Moreover, we also have 
\begin{equation*}
\left\vert \int_{E}\left\vert \nabla u\right\vert ^{2}\langle \nabla
f,\nabla \phi ^{2}\rangle \right\vert \leq C\int_{\left( B_{p}\left(
A+1\right) \backslash B_{p}\left( A\right) \right) \cap E}\left\vert \nabla
u\right\vert ^{2}\rightarrow 0\;\;\text{as}\;\;A\rightarrow \infty .
\end{equation*}%
Here we have used that $u$ has finite Dirichlet integral on the end $E.$ It
is not difficult to see by the co-area formula and since $\left\vert \nabla
u\right\vert $ is bounded we also have 
\begin{equation*}
\int_{M}\left\vert \nabla u\right\vert ^{2}\left\vert \nabla \phi
\right\vert ^{2}\rightarrow 0\;\;\text{as}\;\;T,A\rightarrow \infty .
\end{equation*}%
Hence, letting $A,T\rightarrow \infty $ in (\ref{inequality1}) we get that $%
\left\vert \nabla \left\vert \nabla u\right\vert \right\vert =R\left\vert
\nabla u\right\vert ^{2}=0$. This implies $|\nabla u|=C$. Since the energy
of $u$ on the nonparabolic end is infinite, this implies $u=\mathrm{const}$.
\end{proof}

We now discuss the case when $M$ is parabolic.

\begin{prop}
\label{prop-sec-step} Let $M$ be a parabolic gradient steady K\"{a}%
hler-Ricci soliton with Ricci curvature bounded below and such that for any $%
x\in M$, $Vol\left( B_{x}\left( 1\right) \right) \geq C>0,$ for a uniform $%
C>0$. Then either it is connected at infinity or it splits isometrically as $%
\mathbb{R}\times N,$ for a compact Ricci flat manifold $N$.
\end{prop}

\begin{proof}[Proof of Proposition \protect\ref{prop-sec-step}]
Suppose $M$ has at least two parabolic ends (by assumption all its ends are
parabolic). Let $E$ be one end, then $F:=M\backslash E$ is another end of $%
M. $

There exists a harmonic function $u$ on $M$ such that 
\begin{equation*}
\lim_{x\rightarrow \infty ,\;x\in E}u\left( x\right) =\infty \;\;\text{and\
\ }\lim_{x\rightarrow \infty ,\;x\in F}u\left( x\right) =-\infty .
\end{equation*}%
Applying Theorem 2.1 in \cite{NR} on $E$ and $F$ separately, it follows that 
$\ u$ has bounded gradient on each of the two ends, therefore 
\begin{equation*}
\left\vert \nabla u\right\vert \leq C\ \ \ \text{ on \ }M.
\end{equation*}
Furthermore, with a similar argument as in the proof of (\ref{i1}) in
Proposition \ref{prop-first-step} we obtain that 
\begin{equation*}
\int_{B_{p}\left( r\right) }\left\vert \nabla u\right\vert ^{2}\leq Cr,
\end{equation*}
for any $r>0$ large enough. Applying again Lemma 3.1 in \cite{L} we get that 
$u$ is pluriharmonic i.e. $u_{\alpha \bar{\beta}}=0$. Denoting 
\begin{equation*}
F:=\left\langle \nabla f,\nabla u\right\rangle ,
\end{equation*}
by the argument in Theorem \ref{thm-no-har0}, we get 
\begin{equation*}
\int_{M}|\nabla F|^{2}\phi ^{2}\leq \frac{C}{r^{2}}\cdot
\int_{B_{p}(2r)\backslash B_{p}(r)}|\nabla u|^{2}|\nabla f|^{2}\leq \frac{C}{%
r}\rightarrow 0\,\,\,\mbox{as}\,\,\,r\rightarrow \infty .
\end{equation*}%
This implies that $F$ is constant on $M$ i.e., 
\begin{equation}
\langle \nabla u,\nabla f\rangle =a\in \mathbb{R}.  \label{1}
\end{equation}%
We use again the inequality (\ref{inequality}) proved in Theorem \ref%
{th-steady}, which also holds true in our setting, because $u$ is harmonic.
We now choose the following cut-off $\phi $, defined on the level sets of $u$
(which are compact). For $T$ large enough let 
\begin{equation*}
\phi \left( x\right) =\left\{ 
\begin{array}{c}
0 \\ 
\frac{2T-u}{T} \\ 
1 \\ 
\frac{u+2T}{T} \\ 
0%
\end{array}%
\right. 
\begin{array}{l}
\text{on} \\ 
\text{on} \\ 
\text{on} \\ 
\text{on} \\ 
\text{on}%
\end{array}%
\begin{array}{l}
u\geq 2T, \\ 
T<u<2T, \\ 
-T\leq u\leq T, \\ 
-2T<u<-T, \\ 
u\leq -2T.%
\end{array}%
\end{equation*}%
Observe now that for any $T$ we have 
\begin{gather*}
2\int_{M}\langle \nabla f,\nabla u\rangle \cdot \langle \nabla u,\nabla \phi
^{2}\rangle =4a\int_{M}\phi \langle \nabla u,\nabla \phi \rangle \\
=-\frac{4a}{T^{2}}\int_{T<u<2T}(2T-u)\left\vert \nabla u\right\vert ^{2}+%
\frac{4a}{T^{2}}\int_{-2T<u<-T}\left( u+2T\right) \left\vert \nabla
u\right\vert ^{2} \\
=-\frac{4a}{T^{2}}\left( \int_{T}^{2T}\left( 2T-t\right) dt\right)
\int_{u=t}\left\vert \nabla u\right\vert +\frac{4a}{T^{2}}\left(
\int_{-2T}^{-T}\left( t+2T\right) dt\right) \int_{u=t}\left\vert \nabla
u\right\vert \\
=-2a\int_{u=t}\left\vert \nabla u\right\vert +2a\int_{u=t}\left\vert \nabla
u\right\vert =0.
\end{gather*}%
We have used the co-area formula and, since $u$ is harmonic, $%
\int_{u=t}\left\vert \nabla u\right\vert $ is finite and independent of $t$
for all $t\in \mathbb{R}.$ Moreover, similarly 
\begin{gather*}
\int_{M}\left\vert \nabla u\right\vert ^{2}\langle \nabla f,\nabla \phi
^{2}\rangle =2\int_{M}\left\vert \nabla u\right\vert ^{2}\langle \nabla
f,\nabla u\rangle \phi \phi ^{\prime } \\
=-\frac{2a}{T^{2}}\int_{T<u<2T}\left\vert \nabla u\right\vert ^{2}\left(
2T-u\right) +\frac{2a}{T^{2}}\int_{-2T<u<-T}\left\vert \nabla u\right\vert
^{2}\left( u+2T\right) \\
=-a\int_{u=t}\left\vert \nabla u\right\vert +a\int_{u=t}\left\vert \nabla
u\right\vert =0.
\end{gather*}%
This shows that the last two terms involving $\phi $ in (\ref{inequality})
are in fact zero for all $T.$ Moreover, by co-area formula and since $%
\left\vert \nabla u\right\vert $ is bounded we see that 
\begin{equation*}
\int_{M}\left\vert \nabla u\right\vert ^{2}\left\vert \nabla \phi
\right\vert ^{2}\rightarrow 0\;\;\text{as}\;\;T\rightarrow \infty
\end{equation*}%
Then letting $T\rightarrow \infty $ in (\ref{inequality}) we conclude that $%
\left\vert \nabla \left\vert \nabla u\right\vert \right\vert =R\left\vert
\nabla u\right\vert ^{2}=0$. This implies $|\nabla u|=C$. If $C=0$ this
means $u=\mathrm{const}$ and we are done. Otherwise, $R\equiv 0$, therefore $%
Ric\equiv 0$. This gives now that $f_{ij}=0$, which is known to imply that $%
M $ is isometric to $\mathbb{R}\times N$, for a compact Ricci flat manifold $%
N$. This concludes the proof.
\end{proof}

From Proposition \ref{prop-first-step} and Proposition \ref{prop-sec-step}
we obtain the following.

\begin{thm}
If $M$ is a steady gradient K\"{a}hler-Ricci soliton with Ricci curvature
bounded below and such that for any $x\in M$, $Vol\left( B_{x}\left(
1\right) \right) \geq C>0$ for some constant $C$ independent of $x$, then
either it is connected at infinity or it splits isometrically as $M=\mathbb{R%
}\times N$, for a compact Ricci flat manifold $N$.
\end{thm}

\section{Volume of steady Ricci solitons\label{vol_steady}}

We now study the volume of gradient steady Ricci solitons. Our motivation is
the results obtained in \cite{C, CD, CN} for shrinking Ricci solitons. We
recall that $\left( M,g\right) $ is a complete noncompact gradient steady
Ricci soliton, i.e. 
\begin{equation*}
Ric=Hess\left( f\right) .
\end{equation*}%
It is known that there exists a positive constant $\lambda >0$ such that 
\begin{equation*}
\left\vert \nabla f\right\vert ^{2}+R=\lambda .
\end{equation*}

\begin{thm}
\label{volume} If $(M,g)$ is a gradient steady Ricci soliton there exist
uniform constants $a,c$ and $r_{0}$ so that for any $r>r_{0}$ 
\begin{equation}
ce^{a\sqrt{r}}\geq Vol(B_{p}(r))\geq c^{-1}r.  \label{eq-vol-gr}
\end{equation}
\end{thm}

\begin{proof}[Proof of Theorem \protect\ref{volume}]
We first establish the volume lower bound. If for every $r_{0}$ we have $%
\int_{B_{p}(r_{0})}R=0$ then $R\equiv 0$ on $M$. By 
\begin{equation*}
\Delta _{f}R=R-2\left\vert Ric\right\vert ^{2}
\end{equation*}
this implies $Ric\equiv 0$. In this case it is known that we have (\ref%
{eq-vol-gr}).

Assume there is an $r_{0}>0$ so that $C_{0}:=\int_{B_{p}(r_{0})}R>0$. Then,
since $R\geq 0$, for $r\geq r_{0}$ we have 
\begin{equation*}
C_{0}\leq \int_{B_{p}(r)}R=\int_{B_{p}(r)}\Delta f=\int_{\partial B_{p}(r)}%
\frac{\partial f}{\partial n}\leq \int_{\partial B_{p}(r)}|\nabla f|\leq 
\sqrt{\lambda }\cdot A(\partial (B_{p}(r)),
\end{equation*}%
where we have used $|\nabla f|\leq \sqrt{\lambda }$ on a gradient steady
soliton. This implies for $r\geq r_{0}$ 
\begin{equation*}
Area(\partial B_{p}(r))\geq c>0,
\end{equation*}%
for a uniform constant $c$. If we integrate the previous inequality over $%
[r_{0},r]$ we obtain for $r\geq 2r_{0}$ 
\begin{equation*}
Vol(B_{p}(r))\geq c(r-r_{0})\geq c_{0}\cdot r.
\end{equation*}

We now prove the volume upper bound.

We denote by $dV|_{\exp _{p}\left( r\xi \right) }=J\left( r,\xi \right)
drd\xi $ the volume form of $M$, where $\xi \in S_{p}M$. We will omit the
dependence on $\xi .$ It is known that along a normal minimizing geodesic
starting from $p$, 
\begin{equation*}
\left( \frac{J^{\prime }}{J}\right) ^{\prime }\left( r\right) +\frac{1}{n-1}%
\left( \frac{J^{\prime }}{J}\right) ^{2}\left( r\right) +Ric\left( \frac{%
\partial }{\partial r},\frac{\partial }{\partial r}\right) \leq 0.
\end{equation*}%
We integrate this from $\ 1$ to $r\geq 1$ and use that 
\begin{equation*}
Ric\left( \frac{\partial }{\partial r},\frac{\partial }{\partial r}\right)
=f^{\prime \prime }\left( r\right)
\end{equation*}
to get 
\begin{equation*}
\frac{J^{\prime }}{J}\left( r\right) +\frac{1}{n-1}\int_{1}^{r}\left( \frac{%
J^{\prime }}{J}\right) ^{2}\left( t\right) dt\leq -f^{\prime }\left(
r\right) +C_{0},
\end{equation*}%
for some constant $C_{0}>0$, independent of $r.$ Let us denote%
\begin{equation*}
u\left( t\right) :=\frac{J^{\prime }}{J}\left( r\right) .
\end{equation*}%
Since $f$ has bounded gradient we get, for any $r\geq 1,$ 
\begin{equation*}
u\left( r\right) +\frac{1}{n-1}\int_{1}^{r}u^{2}\left( t\right) dt\leq C.
\end{equation*}%
Notice that by the Cauchy-Schwarz inequality it follows:%
\begin{equation}
u\left( r\right) +\frac{1}{\left( n-1\right) r}\left( \int_{1}^{r}u\left(
t\right) dt\right) ^{2}\leq C.  \label{v3}
\end{equation}%
We claim that for any $r\geq 1,$ 
\begin{equation}
\int_{1}^{r}u\left( t\right) dt\leq \sqrt{\left( n-1\right) Cr}.  \label{v4}
\end{equation}%
To prove this, define 
\begin{equation*}
v\left( r\right) :=\sqrt{\left( n-1\right) Cr}-\int_{1}^{r}u\left( t\right)
dt.
\end{equation*}%
We prove (\ref{v4}) by showing that 
\begin{equation*}
v\left( r\right) \geq 0\text{ \ for all \ }r\geq 1.
\end{equation*}
Clearly, $v\left( 1\right) >0.$ Assume by contradiction that $v$ is not
positive for all $r\geq 1, $ so let $r_{0}>1$ be the first number for which $%
v=0.$

Since $v\left( r_{0}\right) =0,$ it follows that 
\begin{equation*}
\int_{1}^{r_{0}}u\left( t\right) dt=\sqrt{\left( n-1\right) Cr_{0}}.
\end{equation*}
By (\ref{v3}), this implies 
\begin{equation*}
u\left( r_{0}\right) \leq C-\frac{1}{\left( n-1\right) r_{0}}\left(
\int_{1}^{r_{0}}u\left( t\right) dt\right) ^{2}=C-\frac{1}{\left( n-1\right)
r_{0}}\left( n-1\right) Cr_{0}=0.
\end{equation*}%
Consequently, we obtain 
\begin{equation*}
v^{\prime }\left( r_{0}\right) =\sqrt{\frac{\left( n-1\right) C}{2r_{0}}}%
-u\left( r_{0}\right) >0.
\end{equation*}
This implies the existence of a small enough $\delta >0$ such that $v\left(
r_{0}-\delta \right) <v\left( r_{0}\right) =0,$ which contradicts the choice
of $r_{0}.$

We have proved that (\ref{v4}) is true for any $r\geq 1,$ which by the
definition of $u$ means that 
\begin{equation*}
\log J\left( r\right) -\log J\left( 1\right) \leq \sqrt{\left( n-1\right) Cr}%
.
\end{equation*}%
This proves that for any $r\geq 1$ we have an area bound of the form 
\begin{equation*}
Area\left( B_{p}\left( r\right) \right) \leq Ce^{a\sqrt{r}}.
\end{equation*}%
The Theorem is proved.
\end{proof}

Let us note that the lower bound of the volume is sharp. Indeed, the product
of the cigar soliton $\left( \mathbb{R}^{2},\frac{dx^{2}+dy^{2}}{%
1+x^{2}+y^{2}}\right) $with any compact Ricci flat $n-2$ dimensional
manifold is a nonflat steady Ricci soliton with linear volume growth.
However, we do not know any examples of steady Ricci solitons with faster
than polynomial volume growth. While the volume upper bound in Theorem \ref%
{volume} may not be sharp, it is still a very useful information. In fact,
it is crucial in the following estimate for the potential function $f.$

\begin{cor}
\label{potential} Let $\left( M,g\right) $ be a steady nonflat gradient
Ricci soliton. Then there exist $\lambda >0$ and $c>0$ such that for any $%
r\geq 1$ 
\begin{equation*}
\sqrt{\lambda }+\frac{c}{r}\geq \frac{1}{r}\sup_{\partial B_{p}\left(
r\right) }f\left( x\right) \geq \sqrt{\lambda }-\frac{c}{\sqrt{r}}.
\end{equation*}
\end{cor}

\begin{proof}[Proof of Corollary \protect\ref{potential}]
It is known that there exists a constant $\lambda >0$ such that 
\begin{equation*}
R+\left\vert \nabla f\right\vert ^{2}=\lambda .
\end{equation*}
Since $R\geq 0,$ we have $\left\vert \nabla f\right\vert \leq \sqrt{\lambda }%
,$ which proves the upper bound estimate for $f.$

We now show the lower bound. We check directly that 
\begin{equation*}
\Delta e^{f}=\left( \Delta f+\left\vert \nabla f\right\vert ^{2}\right)
e^{f}=\lambda e^{f}.
\end{equation*}%
Integrating this on $B_{p}\left( r\right) ,$ it follows 
\begin{equation}
\lambda \int_{B_{p}\left( r\right) }e^{f}=\int_{B_{p}\left( r\right) }\Delta
e^{f}=\int_{\partial B_{p}\left( r\right) }\frac{\partial }{\partial r}%
\left( e^{f}\right) \leq \sqrt{\lambda }\int_{\partial B_{p}\left( r\right)
}e^{f}.  \label{x1}
\end{equation}%
In the last inequality, we have used that $\left\vert \frac{\partial f}{%
\partial r}\right\vert \leq \left\vert \nabla f\right\vert \leq \sqrt{%
\lambda }.$ Denoting 
\begin{equation*}
w\left( r\right) :=\int_{B_{p}\left( r\right) }e^{f}
\end{equation*}
it follows from (\ref{x1}) that%
\begin{equation*}
\sqrt{\lambda }w\left( r\right) \leq w^{\prime }\left( r\right) .
\end{equation*}%
We integrate this from $1$ to $r$ to conclude that $w\left( r\right) \geq
ce^{\sqrt{\lambda }r},$ for some $c>0.$ By (\ref{x1}) this means that 
\begin{equation*}
\int_{\partial B_{p}\left( r\right) }e^{f}\geq ce^{\sqrt{\lambda }r},\ \ 
\text{for any }r\geq 1.
\end{equation*}%
Combining with our area estimate from Theorem \ref{volume} we get:%
\begin{equation*}
Ce^{a\sqrt{r}}\left( \sup_{B_{p}\left( r\right) }e^{f}\right) \geq \left(
\sup_{B_{p}\left( r\right) }e^{f}\right) Area\left( \partial B_{p}\left(
r\right) \right) \geq \int_{\partial B_{p}\left( r\right) }e^{f}\geq ce^{%
\sqrt{\lambda }r}.
\end{equation*}%
This implies the lower bound for $f$ and proves the Corollary.
\end{proof}

We want to point out that in contrast to shrinking Ricci solitons, where we
have (\ref{eq-cao}), for steady solitons it is not possible to obtain such
estimates for the potential. This is because if $M$ is a gradient steady
Ricci soliton, then so is $\mathbb{R}\times M,$ where the potential is now
constant on $\mathbb{R}$. Thus, in this case the potential does not grow
linearly on $M,$ and in the general setting the result in Corollary \ref%
{potential} seems to be the best we can say. We should also note that
Corollary \ref{potential} was independently proved by P.\ Wu in \cite{Wu}.

\textbf{Acknowledgements.}

We would like to thank Huai-Dong Cao, Peter Li, Mohan Ramachandran and
William Wylie for their interest.

{\tiny DEPARTMENT OF MATHEMATICS, COLUMBIA UNIVERSITY, NEW YORK, NY 10027}%
\newline
{\small E-mail address: omuntean@math.columbia.edu}

\bigskip

{\tiny DEPARTMENT OF MATHEMATICS, UNIVERSITY OF PENNSYLVANIA, PHILADELPHIA,
PA 19104}\newline
{\small E-mail address: natasas@math.columbia.edu}

\end{document}